\documentclass[11pt]{article}
\usepackage{graphicx}
\usepackage{amsfonts}
\usepackage{amsmath}
\usepackage{amssymb}

\usepackage{xspace}
\usepackage{paralist}
\usepackage{ifpdf}
\usepackage{braket}
\usepackage{tabulary}
\usepackage{bbm}
\usepackage{subfig}
\usepackage{algorithm}
\usepackage{algpseudocode}

\usepackage{bm}

\usepackage{soul,color}

\begin{document}
\title{Using Deep Q-Learning to Control Optimization Hyperparameters}
\author{Samantha Hansen \\ IBM T.J. Watson Research Center}
\date{}
\maketitle
\begin{abstract}{We present a novel definition of the reinforcement learning state, actions and reward function that allows a deep Q-network (DQN) to learn to control an optimization hyperparameter. Using Q-learning with experience replay, we train two DQNs to accept a state representation of an objective function as input and output the expected discounted return of rewards, or q-values, connected to the actions of either adjusting the learning rate or leaving it unchanged.
The two DQNs learn a policy similar to a line search, but differ in the number of allowed actions.  The trained DQNs in combination with a gradient-based update routine form the basis of the Q-gradient descent algorithms. To demonstrate the viability of this framework, we show that the DQN's q-values associated with optimal action converge and that the Q-gradient descent algorithms outperform gradient descent with an Armijo or nonmonotone line search.
Unlike traditional optimization methods, Q-gradient descent can incorporate any objective statistic and by varying the actions we gain insight into the type of learning rate adjustment strategies that are successful for neural network optimization.
}
\end{abstract}

\section{Introduction}
This paper demonstrates how to train a deep Q-network (DQN) to control an optimization hyperparameter. Our goal is to minimize an objective function through gradient-based updates of the form
\begin{equation}
x_{t+1} = x_t - \alpha_t g_t \label{g-update}
\end{equation}
where $\alpha_t$ is the learning rate. At each iterate $x_t$, we extract information about the objective derived from Taylor's theorem and line search methods to form a state feature vector. The state feature vector is the input to a DQN and the output is the expected discounted return of rewards, or q-value, connected to the action of increasing, decreasing, or preserving the learning rate. We present a novel definition of the reinforcement learning  problem that allows us to train two DQNs using Q-learning with experience replay \cite{lin1993reinforcement, watkins1992q} to successfully control the learning rate and learn the q-values associated with the optimal actions.

The motivation for this work is founded on the observation that gradient-based algorithms are effective for neural network optimization, but are highly sensitive to the choice of learning rate \cite{lecun2012efficient}. Using a DQN in combination with a gradient-based optimization routine to iteratively adjust the learning rate eliminates the need for a line search or hyperparameter tuning, and is the concept for the Q-gradient descent algorithm. Although we restrict this paper to deterministic optimization, this framework can extend to the stochastic regime where only gradient estimates are available.

We train two DQNs to minimize a feedforward neural network that performs phone classification
in two separate environments. 
The first environment conforms to an Armijo line search procedure \cite{armijo1966minimization, nocedal2006numerical}, either the learning rate is decreased by a constant factor or an iterate is accepted and the learning rate is reset to an initial value. The second environment differs in that the learning rate can also increase and is never reset. 
The trained DQNs are the input to the Q-gradient descent (Q-GD) versions 1 \& 2, and we test them 
against gradient descent with an Armijo or nonmonotone \cite{grippo1986nonmonotone} line search to show that these new algorithms are able to find better solutions on the original neural network, as well as on a neural network that is doubled in size and with three times the amount of data. We also compare how each algorithm adjusts the learning rate during the course of the optimization procedure in order to extract characteristics that explain Q-GD's superior performance.



The paper is organized as follows: in Section~\ref{background} we review reinforcement learning (RL) theory and in Section~\ref{rl-opt} we define the RL actions, state, and reward function for the purpose of optimization. Section~\ref{sec:training} describes the Q-learning with experience replay procedure used to train the DQNs.
In Section~\ref{sec:experiments}, we test the Q-GD algorithms against gradient descent with an Armijo or nonomontone line search on two neural networks that perform phone classification. Section~\ref{sec:related} reviews relevant literature and finally, in Section~\ref{sec:conclusions} we provide concluding remarks and discuss future areas of research.
\\ \\
\noindent
\emph{Notation:} We use brackets indexed by either location or description to denote accessing an element from a vector. For example,  $[s]_i$ denotes the $i^{th}$ element and $[s]_{\textrm{encoding}}$ denotes the element corresponding to description `encoding' for vector $s$.

\section{Review of Reinforcement Learning} \label{background}
Reinforcement learning is the presiding methodology for training an \emph{agent} to perform a \emph{task} within an \emph{environment}. These tasks are characterized by a clear underlying goal and require the agent to sequentially select an \emph{action} based on the \emph{state} of the environment and the current \emph{policy}. The agent learns by receiving feedback from the environment in the form of a \emph{reward}. 

At each time step $t$, the agent receives a representation of the environment's state $s_t\in\mathcal{S}$ and based on the policy $\pi: \mathcal{S}\rightarrow \mathcal{A}$ chooses an action $a_t\in \mathcal{A}$. The agent receives a reward $r_{t+1}$ for taking action $a_t$ and arriving in state $s_{t+1}$. We assume that the environment is a Markov Decision Process (MDP), i.e. given the current state $s_t$ and action $a_t$, the probability of arriving in next state $s_{t+1}$ and receiving reward $r_{t+1}$ does not depend on any of the previous states or actions.

A successful policy must balance the immediate reward with the agent's overall goal. 
RL achieves this via the action-value function $Q^{\pi}(s,a):(\mathcal{S},\mathcal{A})\rightarrow \mathbb{R}$, which is the discounted expected return of rewards given the state, action, and policy,
\begin{equation}
Q^{\pi}(s,a)  = \mathbb{E}_{\pi} \left[ R_{t+1} \ \Bigr\rvert  \ s_t = s, \ a_t =a \right]
\end{equation}
where
\begin{equation}
R_{t+1} = r_{t+1} + \sum_{k=1}^{T-t-1} \gamma^k r_{t+1+k}, \ 0 < \gamma \leq 1, \label{returnReward}
\end{equation}
$T$ is the maximum number of time steps and the expectation is taken given that the agent is following policy $\pi$.
The optimal action-value function, $Q^{*}(s,a) =\max_{\pi} Q^{\pi}(s,a) $ satisfies the Bellman equation,
 \begin{equation}
Q^{*}(s,a) = \mathbb{E}_{\pi^*}\left[r_{t+1} + \gamma\max_{a'\in\mathcal{A}}Q^{*}(s_{t+1},a') \  \Bigr\rvert \ s_t = s, \ a_t =a \right] \label{bellman}
\end{equation}
which provides a natural update rule for learning.  At each time step the effective estimate $\hat{y}_t$ and target $y_t$ are given by
\begin{equation}
 \hat{y_t} = Q_t(s_t,a_t), \ \ \ y_t = r_{t+1} + \gamma\max_{a'\in\mathcal{A}}Q_t(s_{t+1},a') \label{Qlearning1} 
\end{equation}
 and the update is based on their difference; this method is referred to as Q-learning. Notice that the estimate/target come from LHS/RHS of (\ref{bellman}) and will both continue to change until $Q$ converges. 
 
 For finite number of states and actions $Q$ is a look-up table. When the number of states is too large or even infinite, the table is approximated by a function. In particular, when the action-value function is a neural network it is referred to as a deep Q-network (DQN). A practical choice is to choose a network architecture such that the inputs are the states and the outputs are the expected discounted return of rewards, or q-value, for each action. We only consider the case of using a DQN and henceforth use the notation 
 \begin{equation}
 Q(s;\theta): \mathbb{R}^{|\mathcal{S}|}\rightarrow \mathbb{R}^{|\mathcal{A}|} \label{dqn}
 \end{equation}
to denote that the DQN is parameterized by weights $\theta$.

The weights are updated by minimizing the $\ell_2$ norm between the estimate and target, $\left\| \hat{y}_t - y_t \right\|_2^2$, yielding iterations of the form
\begin{equation}
 \theta \leftarrow  \theta- \beta ( \hat{y}_t- y_t)\nabla_{\theta} Q(s_t; \theta)\label{Qlearning2}
\end{equation}
where $\beta$ is the learning rate and
\begin{equation}\label{q-update}
\hat{y}_t =  Q(s_t; \theta), \qquad [y_t]_a = \begin{cases}
r_{t+1} + \mathbbm{1}_{t\neq T-1} \left(\gamma\max_{a'\in\mathcal{A}}[Q(s_{t+1};\theta)]_{a'} \right)& \ \  a=a_t \\
[Q(s_t; \theta)]_a & \ \ a\neq a_t.
\end{cases}\end{equation}
For the last action, only the reward is present in the target definition and for the non-chosen actions, the targets are set to force the error to be zero.

The action at each time step is chosen based on the principle of exploration versus exploitation. Exploitation takes advantage of the information already garnered by the DQN while exploration encourages random actions to be taken in prospect of finding a better policy. We employ an $\epsilon$-greedy policy which chooses the optimal action w.r.t the DQN's q-values with probability $1-\epsilon$ and randomly otherwise:
\begin{equation} \label{eps-greedy}
a_t =  \begin{cases} \arg\max_a [Q(s_t;\theta)]_a & r \geq \epsilon \\ \textrm{randomly chosen action} & r < \epsilon \end{cases}
\end{equation}
where $r\sim U[0,1]$. 

Equation (\ref{eps-greedy}) is the effective policy since it maps states to actions. Q-learning is an off-policy procedure because it follows a non-optimal policy (with probability $\epsilon$ a random action is taken) yet makes updates to the optimal policy, as illustrated by the $\max$ term in (\ref{q-update}). For a comprehensive introduction to RL, see \cite{sutton1998reinforcement}.
\section{Reinforcement Learning for Optimization} \label{rl-opt}

In this section, we outline the environment, state, actions, and reward function that define the reinforcement learning problem for the purpose of optimization.


 %


\subsection{Actions}
We present two procedures for adjusting the learning rate and show how they are implemented in practice. The first strategy mimics an Armijo line search \cite{armijo1966minimization, nocedal2006numerical} in that the learning rate is reset to an initial value after accepting an iterate and can only henceforth be decreased. The second strategy permits the
learning rate to increase or decrease and is never reset. The two methods are outlined in Algorithm~\ref{alg:a2} and are referred to as Q-gradient descent (Q-GD) versions 1 \& 2, respectively.

Q-GD is a gradient descent optimization procedure that uses a trained DQN to determine the learning rate. The Q-GD inputs are an initial iterate and learning rate  $x_1$ and $\alpha_c$, trained DQN $Q(s;\theta)$, and maximum number of time steps $T$.  We use the notation $x_t$ to denote the candidate iterate, which changes at every time step, and $\bar{x}$ to represent an accepted iterate with associated decent direction $d(\bar{x})$.
In steps \ref{a-s1} and \ref{a-s2}, a state feature vector representative of the objective (discussed in the next section) is formed and passed through the DQN to determine the action. After the action is taken, the candidate iterate is updated in step~\ref{s-update}. 

When a good initial learning rate is known then the first version is preferable, e.g. $d(\bar{x})$ is the Newton direction
and $\alpha_c = 1$ for convex $f$. For non scale-invariant search directions, such as the gradient direction, the second version is advantageous. 

\begin{algorithm}[!]
\caption{Q-gradient descent versions 1 \& 2}
 \label{alg:a2}
\textbf{Input:} initial iterate $x_1$, initial learning rate $\alpha_{c}$, trained DQN $Q(s; \theta)$, number of time steps $T$
\begin{algorithmic}[1]
\State Set $\bar{x} = x_1$, $d(\bar{x}) = -\nabla f(x_1)$, $\alpha_1 =  \alpha_{c}$
\For{$t=1, \ldots, T$} 
\State Compute state feature vector $s_t$ \label{a-s1}
\State $a_t = \arg\max_a [Q(s_t;\theta)]_a$ \label{a-s2}
\If{$ a_t = a_{\textrm{half}}$ }
\State $\alpha_{t+1} =  \frac{1}{2}\alpha_t$ 
\ElsIf{$ a_t = a_{\textrm{double}} $} \Comment{Only for version 2}
\State $\alpha_{t+1} =  2 \alpha_t$ 
\ElsIf{$a_t = a_{\textrm{accept}}$}
\State $\bar{x} = x_{t}, \  d(\bar{x}) = -\nabla f(\bar{x}), \ \alpha_{t+1} = \begin{cases} \alpha_{c} \ & \textrm{version 1} \\ \alpha_t \ & \textrm{version 2} \end{cases} $ \Comment{Update accepted iterate} \label{a-sa}
\EndIf
\State $x_{t+1} =\bar{x} + \alpha_{t+1} d (\bar{x})$ \label{s-update} \Comment{Update candidate iterate}
\EndFor
    \State \Return $x^* = x_T$         
  \end{algorithmic}
\end{algorithm}

\subsection{Environment and State} \label{sec:EandS}
The environment is a combination of the objective function $f: \ \mathbb{R}^n\rightarrow\mathbb{R}$ and set of allowed actions and
needs to be formulated as a MDP in order for the Q-learning algorithm to operate. The Markov condition could be satisfied by including the initial iterate, and the current, as well as all proceeding learning rates and descent directions into the state definition. However, for objective functions with large number of variables such an approach is computationally prohibitive and would severely limit the trained DQN's ability to generalize to a broader family of functions.
We seek to define the state
such that it characterizes the objective function at a given iterate, contains some history, and is universal to all functions.
We use a nonmontone line search as a starting point since it provides an effective criteria for determining the learning rate that is independent of function variable size or type.


A nonmonotone line search chooses the learning rate such that the new iterate is sufficiently less than the maximum objective value of the past $M$ iterates,
\begin{equation}
f(x_{t} + \alpha_td_t) \leq \max_{i=t, \ldots, t-M+1} f(x_i) +c\alpha_t d_t^T\nabla f(x_t) , \ \ c > 0. \label{ls_nm}
\end{equation}
This suggests that the state features needed in order to determine the learning rate are the current learning rate, candidate iterate objective value, max objective from the past $M$ steps, and the dot product between the descent direction $d_t$ and gradient $\nabla f(x_t) $. Although this feature set would neither satisfy the Markov property nor completely capture the objective, updates based on (\ref{ls_nm}) work well in practice and we use these statistics as motivation for the state features.

We employ an encoding that indicates whether the candidate iterate is higher/lower than the M lowest achieved objective values.
Let $F^{t-1}_{M}$  be a list of the $M$ lowest objective values obtained up to time $t-1$, the state encoding is given by
\begin{equation} \label{encoding}
[s_t]_{\textrm{encoding}} = \begin{cases}
1 & \ \ f(x_{t}) \leq \min(F^{t-1}_{M}) \\
0 & \ \  \min(F^{t-1}_{M})  < f(x_{t}) \leq \max(F^{t-1}_{M}) \\ 
-1 & \ \ \textrm{otherwise}.
\end{cases}\end{equation} 
The number of function evaluations must also be a state feature since the states wouldn't otherwise be stationary and the maximum number of time steps $T$ designates an absorbing state. Based on RPROP \cite{riedmiller1993direct}, the final state feature is a measure of alignment between successive descent directions
\begin{equation}
[s_t]_{\textrm{alignment}} = \frac{1}{n}\sum_{i=1}^n \textrm{sign}([d_t]_i [d_{t-1}]_i ). \label{algn}
\end{equation}
In summary there are six features: current learning rate, objective value, dot product between the search direction and gradient, min/max encoding (\ref{encoding}), number of function evaluations, and alignment measure (\ref{algn}).

For the purpose of making the state features independent of the specific objective function, all of the features are transformed to be in the interval $[-1,1]$. For each feature $[s]_i$, a maximum and minimum value is estimated so that
\begin{equation}
[\hat{s}]_i = 1 - 2([s]_i -[s_{\textrm{min}}]_i)/([s_{\max}]_{i} - [s_{\textrm{min}}]_i). \label{transform-x}
\end{equation}
Additionally, since the objective values and gradient norms both converge towards a lower bound $c_i$, these features are transformed twice. First via $[s]_i \leftarrow 1/([s]_i - c_{i})$ and they by (\ref{transform-x}), where $c_i$ is set to 0 for the gradient norm and an objective lower bound $f_{lb}$ for the function values. In general, $f_{lb}$ can be set to zero for objectives that are a sum of loss functions.
 
\subsection{Reward Function}
\begin{figure}[!]
 \raggedleft
\includegraphics[scale=.45]{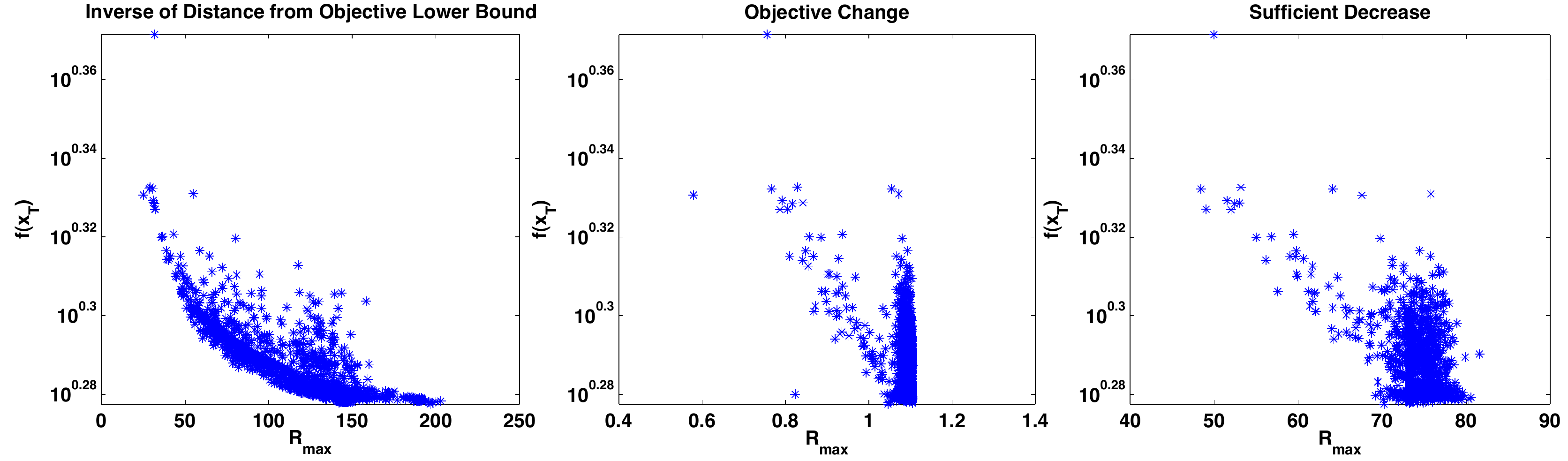}
\caption{Comparison of reward functions. The images plot $f(x_T)$ versus $R_{\max}= \max_t R_t$  for reward functions
defined by $r_{id}$ (inverse distance from objective lower bound),  $r_{oc}$ (objective change) and $r_{sd}$ (sufficient decrease) given by equations \ref{reward} and \ref{otherrewards}.
Only $r_{id}$, shown in the leftmost graph, has the highest $R_{\max}$ values concentrated towards lowest final objective values.}
\label{rewardImage}
\end{figure}

 The reward function is crucial in ensuring that the DQN learns a policy consistent with the goal of finding the lowest objective value in the fewest number of steps, and we define it as the inverse distance from the objective lower bound,
 \begin{equation}
r_{id}(f,x_t) = \frac{c}{f(x_t)-f_{lb}}, \ \ c > 0, \ \ f_{lb} < f(x) \ \  \forall x. \label{reward}
 \end{equation}
The reward function (\ref{reward}) is strictly positive and asymptotes as $f$ approaches
the lower bound. 

We tested reward functions based on a sufficient decrease condition or change in objective value between successive iterates,
\begin{equation}
r_{sd}(f,x_t) = 1_{ f(x_{t-1}) \geq 1.001f(x_t)}, \qquad  r_{\textrm{oc}} (f,x_t) = f(x_{t-1}) - f(x_t) \label{otherrewards}
\end{equation}
and found that they did not adequately capture the optimization goal. To compare the different reward functions we plotted $f(x_T)$ 
against $R_{\max} = \max_t R_t$; for each training episode of DQN v1 we recorded the sequence of objective values ($f(x_T)$ being the objective value at the last time step) and used this information to calculate $R_{\max}$ for each reward function. Figure~\ref{rewardImage} shows that reward functions based on sufficient decrease or objective change 
yield high $R_{\max}$ values for suboptimal final solutions. The main difference between the reward functions is that (\ref{reward}) is based on degree of difficulty in decreasing the objective and will generate the highest rewards during the final time steps.

\section{Training}\label{sec:training}
This section outlines the Q-learning with experience replay method used to train DQN versions 1 \& 2  \cite{watkins1992q, lin1993reinforcement}. 
Algorithm~\ref{alg:Q-learning} exhibits the overall procedure, but omits some of the specific details, which are discussed in the subsections for the sake of clarity.
Note that updates  w.r.t. $f(x)$ are explicitly shown and are indexed by the time step $t$ while the DQN update in step~\ref{dqn-update}  is referenced via equation (\ref{minibatch-update}) and is implicitly indexed by the time step and episode. 

The DQN learns how to minimize the function $f(x)$ through repeated attempts, called learning episodes. For each learning episode, the $x$ iterate is set to an initial value and the DQN then has $T$ time steps to find the lowest objective value. An alternative approach for limiting the number of time steps is to end the episode once the objective has decreased past a certain threshold. Both approaches force the DQN to learn a trade off between finding a good learning rate and exploring the space. Restricting the number of time steps reflects real world applications where there are computational and time constraints and also does not require a-priori knowledge of the objective function.


\subsection{Experience Replay}
An experience consists of a $(s_i, a_i, r_{i+1}, s_{i+1})^j$ tuple for some episode $j\in [1, e]$ at time step $i\in[M-1,T]$, where $M$ and $e$ are in Algorithm~\ref{alg:Q-learning} steps \ref{epi} and \ref{fortime}. These tuples are stored in a memory of experiences $\mathcal{E}$. Instead of updating the DQN with only the most recent experience, a subset $\mathcal{S}\subset\mathcal{E}$ of experiences are drawn from memory and used as a mini-batch to update the DQN:
\begin{equation}
 \theta \leftarrow  \theta- \frac{\beta}{|\mathcal{S}|}\sum_ {(s_i, a_i, r_{i+1}, s_{i+1})^j\in\mathcal{S}}( \hat{y}_i- y_i)\nabla_{\theta} Q(s_i; \theta) \label{minibatch-update}
 \end{equation}
where the estimate $\hat{y}_i$ and target $y_i$ are given via (\ref{q-update}).

The $A$ most recent episodes along with the top $B$ best games (in terms of $R_{\max}$ value) are stored in memory. At each DQN update (step~\ref{dqn-update}) the subsample $\mathcal{S}$ is formed by randomly drawing experiences  from $\mathcal{E}$ and an experience from each of the top $B$ best games. Adding randomly drawn experiences to the mini-batch helps prevent the DQN from over learning during a particular time and episode.


\subsection{Training Specifications}\label{sec:ts}
The Q-learning input parameters in Algorithm~\ref{alg:Q-learning}  for both DQN versions 1 \& 2 were fixed as follows: the discount factor was set to $\gamma = .99$ and the exploration probability $\epsilon$ was initially set to 1 then uniformly decayed to $.1$ over the first 100 episodes. For experience replay, $A=45$, $B=5$, and the mini-batch size was set to $|\mathcal{S}| = 32$. Additionally, for the first 50 episodes the top $B$ best games were not used in the mini-batch sample. The constants $c_1$ and $c_2$ used to calculate the reward (see steps \ref{r1} and \ref{r2}) were fixed as .1 and .12, respectively.  The total number of episodes $E$ is 150K for version 1 and 400K for version 2.

The objective input parameters in Algorithm~\ref{alg:Q-learning} consist of the objective function $f(x)$ with lower bound $f_{lb}$, initial weights $x_1$, initial learning rate $\alpha_c$, encoding memory $M$, and the total number of time steps $T$. The objective function has the form
\begin{equation}
\frac{1}{N} \sum_{i=1}^N \ell(h(z_i; x),t_i) \label{f-form}
\end{equation}
where $z_i$ is an acoustic feature vector with phonetic label $t_i$, $\ell(\cdot)$ is a cross entropy loss, and $h(z; x)$ is a feedforward neural network parameterized by $x$ with sigmoid activations and a softmax function at the output layer. We set the input objective function to $f_{train}$, which has a neural network architecture $65\times16\times8\times42$ and $N=5000$ data points. The number of time steps is $T=1000$ and $M=3$. 


At the start of each episode, the $x$ iterate is reset to $x_1$ and is updated for the first $M$ time steps using the initial learning rate (step~\ref{firstM}) in order to form the first state feature vector.
In steps \ref{getState} and \ref{getAction}, the six state features form the input to the DQN and the resulting action is determined by an $\epsilon$-greedy policy.
Based on the action, the learning rate is either modified, step \ref{step-half} or \ref{step-double},  or the current iterate is accepted and a new gradient direction is calculated, step~\ref{step-accept}. The iterate $x_t$ is updated in step~\ref{xupdate} and this causes the environment to change to the next state (step~\ref{newstate}). The reward for arriving to state $s_{t+1}$ is calculated using either the objective value at the new iterate (step~\ref{r1}) or the previous iterate (step~\ref{r2}) for when the action is to accept. As an aside, 
we found it beneficial to calculate the reward for each action at the last time step since the
targets associated with absorbing states do not change during training and thus play a vital role for propagating back information.
The tuple $(s_t, a_t, r_{t+1},s_{t+1})^e$ forms an experience and is added to memory $\mathcal{E}$ (step~\ref{addexp}). In addition to the current experience, a random subset of experiences are drawn and used to form a mini-batch update for the DQN (step~\ref{dqn-update}).


Special modifications were needed for training DQN v2 since one of its actions permits the learning rate to increase. Too large of a learning rate resulted in updates that caused the objective function to diverge and consequently produce state vectors with infinite features. To prevent this from happening, we used a maximum and minimum learning rate as part of the training procedure. If DQN v2 attempted to increase/decrease the learning rate above/below these values then it would receive a reward of -1 and the episode would terminate early. In addition, we employed an rmsprop update procedure for training DQN v2 \cite{tieleman2012lecture}.

DQN versions 1 \& 2 have an architecture of $6\times 32 \times 16\times |\mathcal{A}|$ with sigmoid activations for the hidden layers and an identify activation for the last layer.
The initial learning rate was set to $\alpha_c=4$ for version 1 and $\alpha_c=2$ for version 2. Additionally, for version 2 only learning rates in the range $[.01, 8]$ were allowed.

\begin{algorithm}
  \caption{Q-Learning with Experience Replay}
  \label{alg:Q-learning}
  \raggedright
  \textbf{Objective Parameters:} $f$, $f_{lb}$, $x_1$, $\alpha_c$, $M$, $T$\\
  \textbf{Q-Learning Parameters:} E, $\theta_0$, $\gamma$, $\epsilon$, $c_1$, $c_2$, $\beta$ \\
  \begin{algorithmic}[1]
  \State $\theta \leftarrow \theta_0$
    \For{$e =1, \ldots, E$}    \Comment{For each learning each episode}   \label{epi}
      \For{$t=1 \ldots, M-1$} \label{fortime} 
      \State $x_{t+1} = x_t - \alpha_c \nabla f(x_t)$ \label{firstM}
      \EndFor
       \State set $\bar{x} = x_{M}$, $d(\bar{x}) = -\nabla f(x_{M})$, $\alpha_{M} =  \alpha_{c}$
      \For{$t=M, \ldots, T$}
      \State Generate state feature vector $s_t$ \label{getState}
      \State Choose action $a_t$ according to $\epsilon$-greedy policy (\ref{eps-greedy}) \label{getAction}
            \If{$ a_t = a_{\textrm{half}}$} \label{step-half}
          \State $\alpha_{t+1} =  \frac{1}{2}\alpha_t$ 
     \ElsIf{$ a_t = a_{\textrm{double}} $} \Comment{Only for version 2} \label{step-double}
          \State $\alpha_{t+1} =  2 \alpha_t$ 
    \ElsIf{$a_t = a_{\textrm{accept}}$}
         \State $\bar{x} = x_{t}, \  d(\bar{x}) = -\nabla f(\bar{x}), \   \alpha_{t+1} = \begin{cases} \alpha_{c} \ & \textrm{version 1} \\ \alpha_t \ & \textrm{version 2} \end{cases}$ \label{step-accept}
    \EndIf
    \State $x_{t+1} =\bar{x} + \alpha_{t+1} d (\bar{x})$ \label{xupdate}
    \State  Generate state feature vector $s_{t+1}$ \label{newstate}
    \If{$a_t \neq a_{\textrm{accept}}$}       
          \State $r_{t+1} = c_1 / (f(x_{t+1}) - f_{lb})$ \label{r1}
     \ElsIf{$a=a_{\textrm{accept}}$}
         \State $r_{t+1} = c_2/ (f(\bar{x}) - f_{lb})$ \label{r2}
      \EndIf
      \State Add experience $(s_t, a_t, r_{t+1}, s_{t+1})^e$ to memory $\mathcal{E}$ \label{addexp}
      \State Sample $\mathcal{S} \in\mathcal{E}$ and update $\theta$ via (\ref{minibatch-update}) \label{dqn-update}
    \EndFor
         \EndFor
    \State \Return $\theta$         
  \end{algorithmic}
\end{algorithm}

\section{Experiments} \label{sec:experiments}
The trained DQNs along with the initial learning rates $\alpha_c$ are the input to the Q-gradient descent algorithms versions 1 \& 2 outlined in Algorithm~\ref{alg:a2}. Since there are no theoretical guarantees that the DQNs would find a good policy or converge, we demonstrate that Q-GD versions 1 \& 2 are effective algorithms by comparing them against gradient descent with an Armijo or nonmonotone line search and show that the DQN q-values associated with the optimal actions converge to the discounted return of rewards at each time step.

The line search algorithms operate under the same rules as Q-GD v1, but an iterate is accepted only if (\ref{ls_nm}) is satisfied. We set  $c=10^{-4}$ and $M=3$ for nonmonotone and, by definition, $M=1$ for Armijo.
\subsection{Results on Train Function} \label{sec:experimentstrain}

\begin{figure}[H]
	\centering
	\subfloat[Objective Value versus Time Step]
	{ \includegraphics[scale=.50]{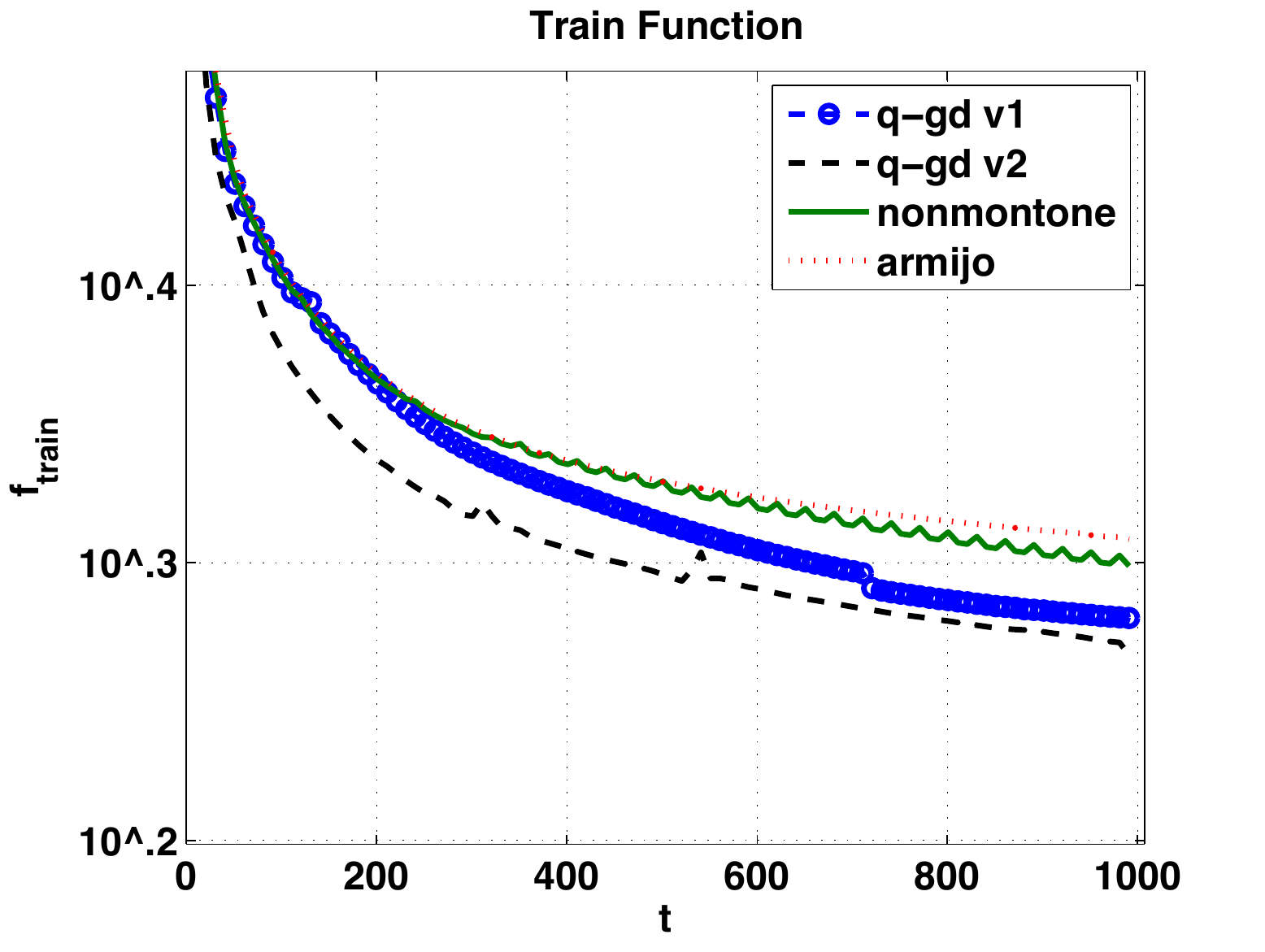}}
	\subfloat[Learning Rate versus Time Step]
	{\hbox{\hspace{-0.95ex}\includegraphics[scale=.50]{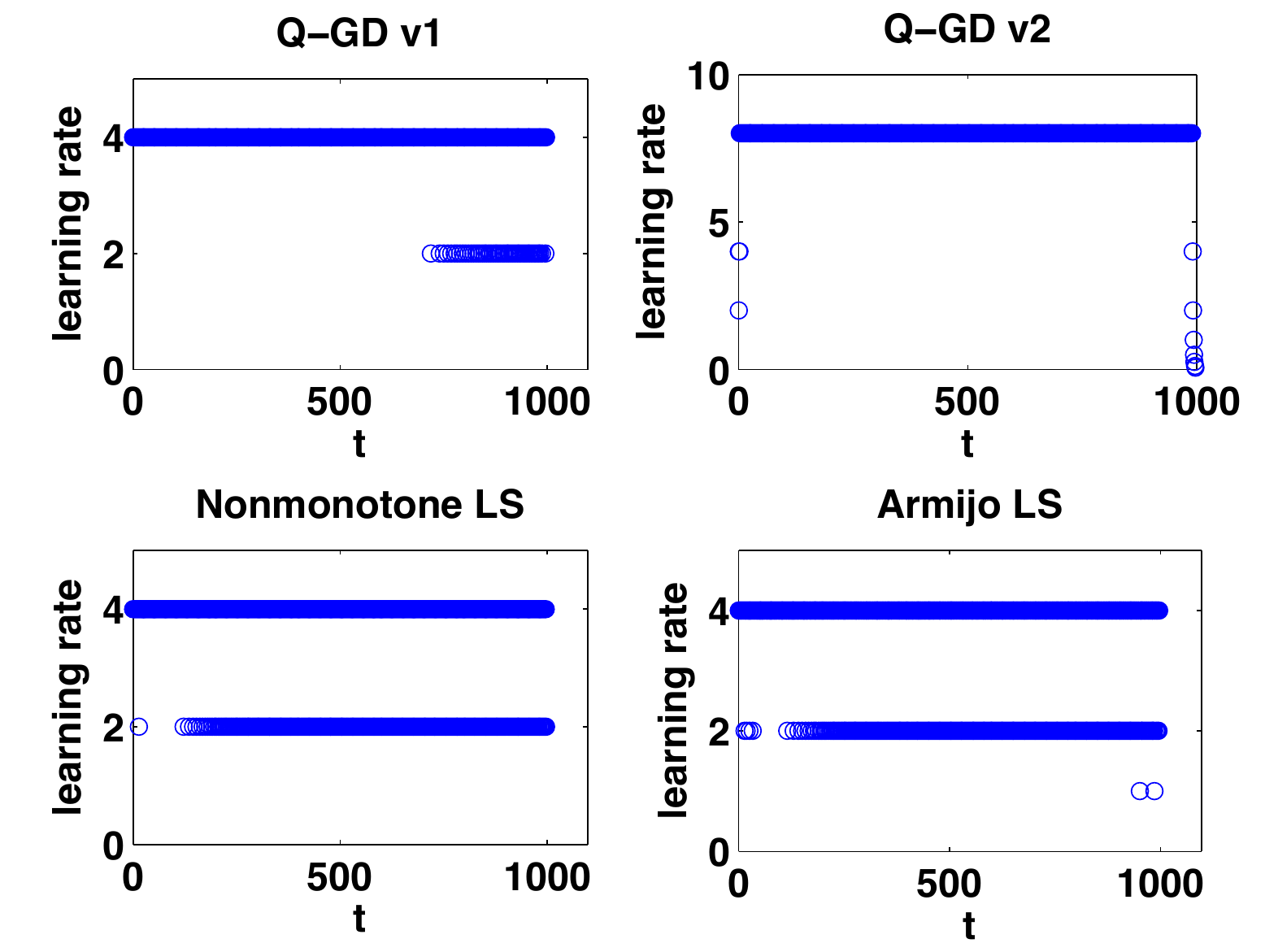}}}
		\caption{Comparison of Q-GD versions 1 \& 2 and gradient descent with a nonmonotone or Armijo line search on train function.}
	\label{fig:ftrain}
\end{figure}

We first compare Q-GD versions 1 \& 2 and gradient descent with an Armijo or nonmonotone line search on the function used to train DQN versions 1 \& 2; $f_{train}$ has the form (\ref{f-form}) with $N=5000$ and feedforward neural network architecture $65\times16\times8\times42$.
Figure~\ref{fig:ftrain}a demonstrates their performance in minimizing $f_{train}$ and figure~\ref{fig:ftrain}b plots the learning rate at each time step.
After 1000 time steps, the final objective values are 1.86, 1.91, 1.98, and 2.04 for Q-GD v2, Q-GD v1, nomonotone, and Armijo, respectively.

The plots of the learning rates illuminate why the Q-GD algorithms are superior. Q-GD v2 has the advantage that it can increase the learning rate and its policy for minimizing the train function was very simple: it increased the learning rate from 2 to 8 during the first initial time steps and then left the learning rate unchanged until decreasing it at each of the last seven time steps. Q-GD v1 offers a fairer comparison to the Armijo and nonomontone line searches since the algorithms all follow the same structure: every time an iterate is accepted the learning rate is reset to 4 and can only then be decreased by a factor of two. The notable difference between Q-GD v1 and the line search algorithms is the frequency in which the learning rate is decreased.
Q-GD v1 decreased the learning rate $5.1\%$ of the time while the Armijo and nonmonotone line searches decreased the learning rate $36.4\%$ and $27.3\%$ of the time. Q-GD v1 also only decreased the learning rate during the final quarter of the optimization procedure.

 The learned policies illustrate that a good initial learning rate is more important than a line search procedure for fast initial objective decrease. Also, it is beneficial to decrease the learning rate more aggressively during the final time steps. Unlike the line searches, the Q-GD algorithms have knowledge of when the optimization procedure is going to end (since the number of time steps is an input parameter) and can act adjust the learning rate accordingly. 


\subsection{Generalization Ability} \label{sec:experimentstest}
\begin{figure}[H]
	\centering
	\subfloat[Objective Value versus Time Step]
	{ \includegraphics[scale=.50]{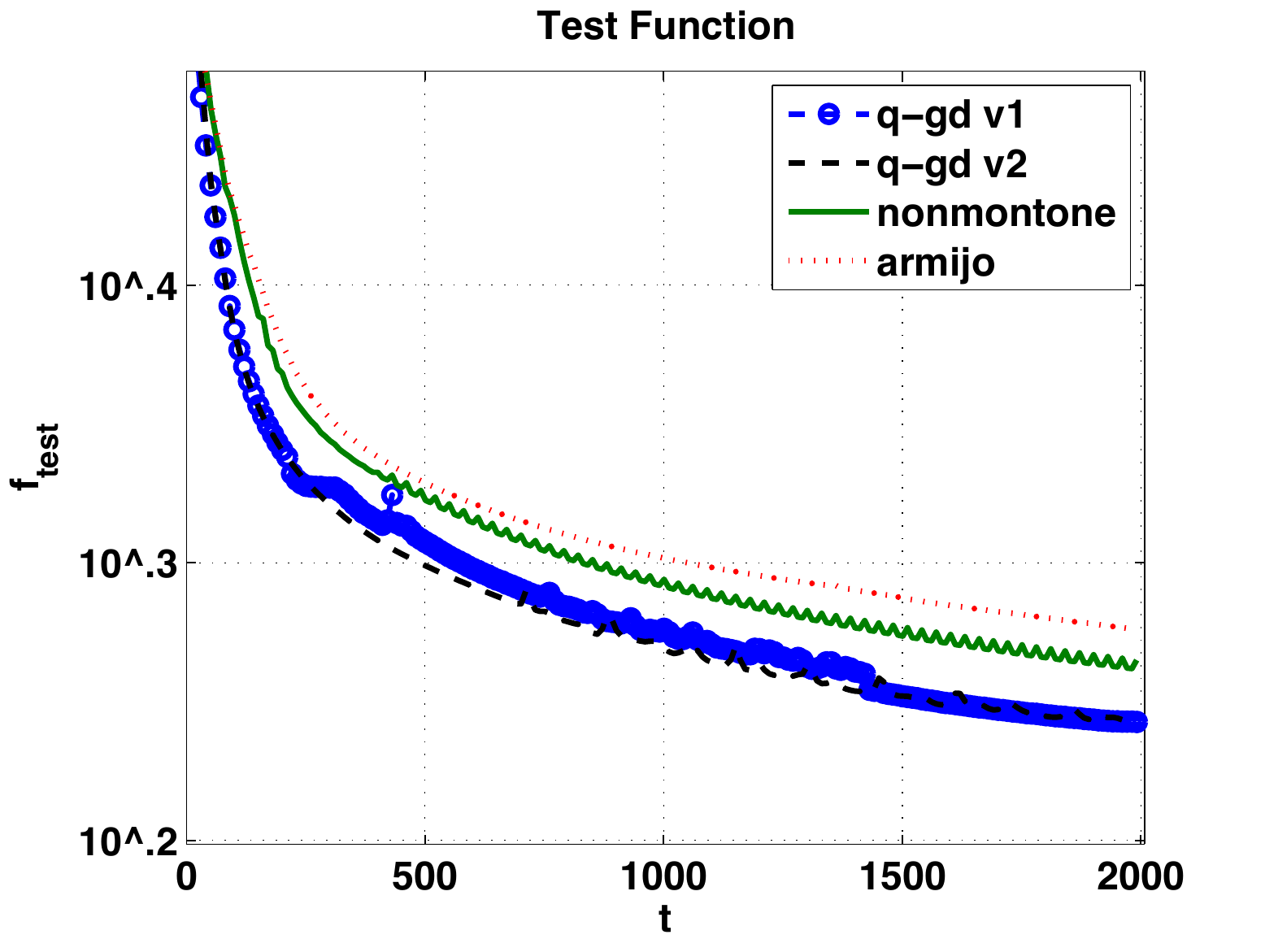}}
	\subfloat[Learning Rate versus Time Step]
	{\hbox{\hspace{-0.95ex}\includegraphics[scale=.50]{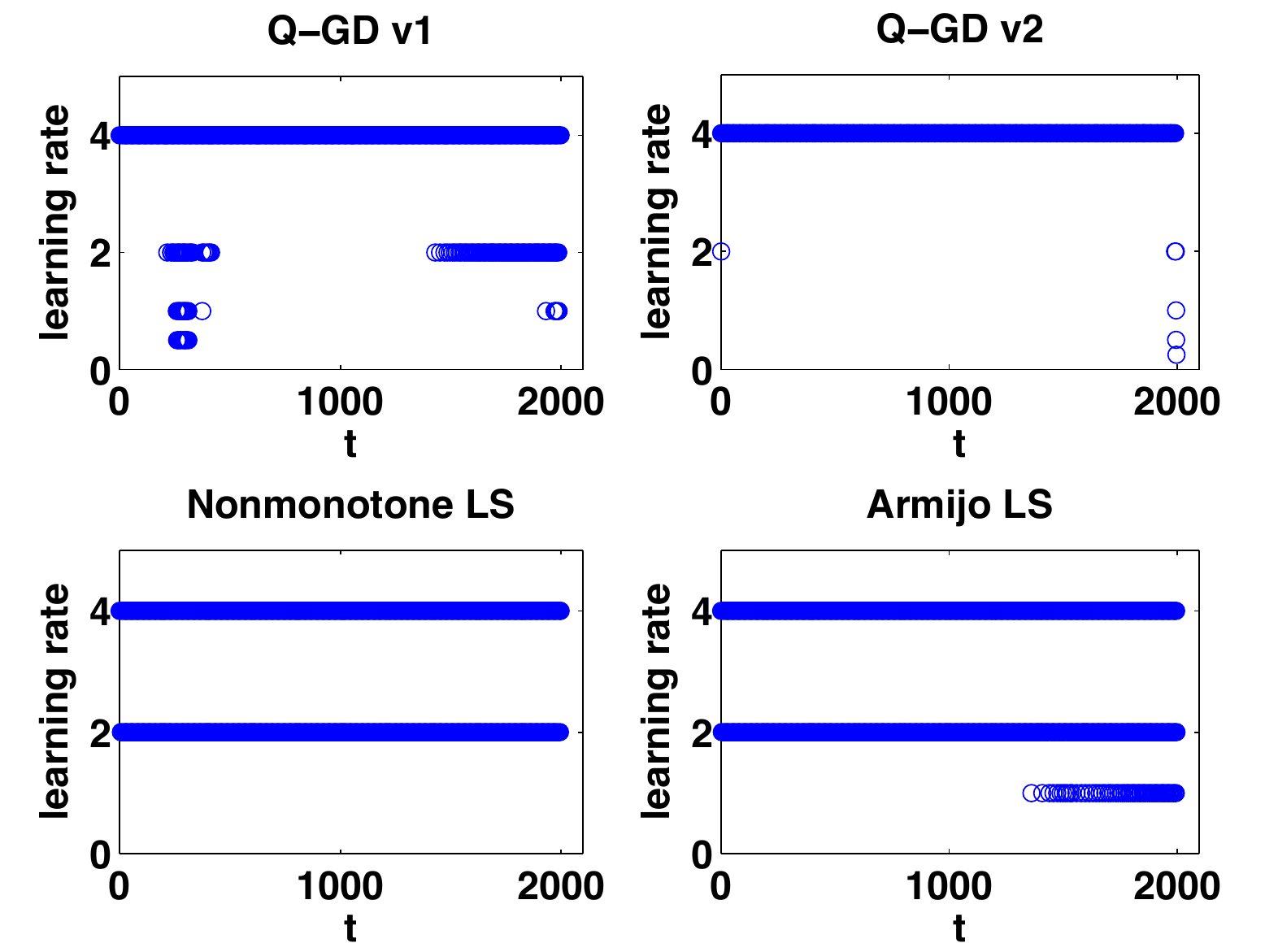}}}
		\caption{Comparison of Q-GD versions 1 \& 2 and gradient descent with a nonmonotone or Armijo line search on test function.}
	\label{fig:ftest}
\end{figure}

We next test to determine if the strategy learned by DQN versions 1 \& 2 on the train function also works for a new, but related function. The test function has the same form as the train function, but with three times the amount of data and double the number of variables (\ref{f-form}) with $N=15000$ and architecture $65\times32\times16\times42$. The purpose of this configuration is to show that we can train the DQN using a small problem and later implement it on larger problems in terms of both variable size and data. We also increased the number of time steps from $1000$ to $2000$.

Figure~\ref{fig:ftest} exhibits how Q-GD versions 1 \& 2 and the nonmonotone and Armijo line search algorithms measure on the test function.
In figure~\ref{fig:ftest}a, we observe that the algorithms retain their relative ordering regarding objective decrease in a fixed number time steps; the final values are 1.73, 1.74, 1.84 and 1.89 for Q-GD v2, Q-GD v1,  nonomonotone and Armijo, respectively. The gap in performance between Q-GD versions 1 \& 2 reduced, showing that Q-GD v1 was more adapt at generalizing to a new function. As with the train function, both Q-GD versions 1 \& 2 decreased the learning rate less frequently than either the nonmonotone or Armijo line searches.
However, both Q-GD versions were more cautious using a higher learning rate at the start of the of the optimization procedure.
Q-GD versions 1 \& 2 maintained their underlying strategies, except version 1 chose to decrease the learning rate during the first quarter and version 2 only initially increased the learning rate to 4 (as opposed to 8).

Overall, these results show that Q-GD versions 1 \& 2 were robust when given a new, larger function and used over a longer number of time steps.


\subsection{Convergence of DQN Q-values}\label{sec:stateFeature}
\begin{figure}[!]
\centering 
\includegraphics[scale=.450]{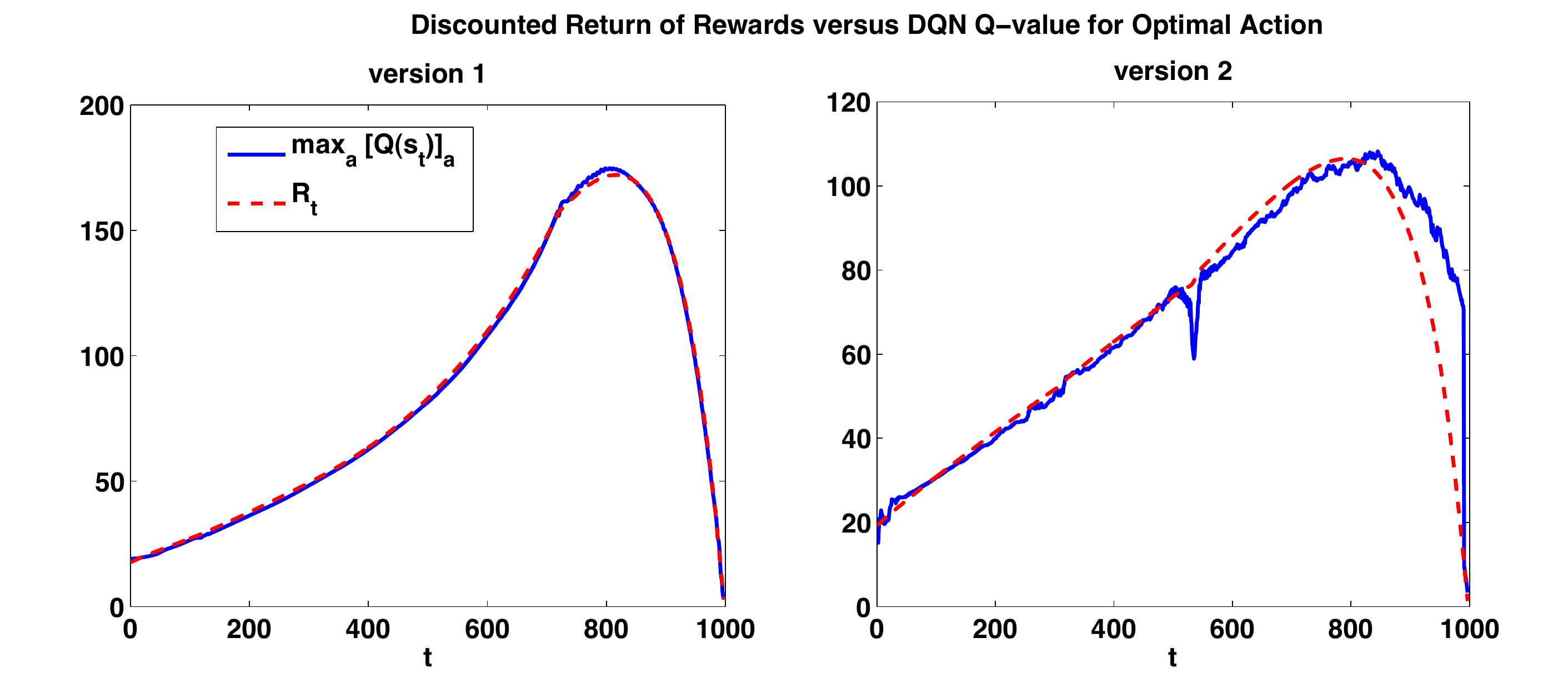}
\caption{Plot of DQN versions 1 \& 2 predicted q-value for optimal action versus the discounted return of rewards (\ref{returnReward}) at each time step $t$ on train function. }
\label{fig:qvalues}
\end{figure}

The purpose of this section is to show that the six state features detailed in Section~\ref{sec:EandS} are rich enough for the DQN to discriminate states in order to learn the q-values associated with the optimal actions. We also demonstrate the effect of individually zeroing out the state features for Q-GD version 1 on the train function.

For the final episode, we recorded the q-value associated with the selected action (no longer using an $\epsilon$-greedy procedure) and resulting reward at each time step in order to compare the DQN predicted q-values against the discounted return of rewards, defined by (\ref{returnReward}). Figure~\ref{fig:qvalues} shows that DQN v1's q-values converged to the discounted return of rewards while DQN v2 found the overall shape of the distribution. Even though DQN v2 was trained with more episodes (400K versus 150K), the addition of one extra action exponentially increases the search space, creating a much more difficult problem.

To investigate how the state features influence the Q-GD algorithms, we
ran Q-GD v1 with either the objective value, gradient norm, or alignment measure set to zero; since the features are transformed to lie in the interval $[-1,1]$ this corresponds to fixing a given feature at its median value. We left the learning rate, objective encoding, and number of time steps unchanged as they are arguably the bare minimum inputs needed to satisfy the Markov property.

Table~\ref{table:features} reports the final objective value and the ratio of halving the learning rate or accepting an iterate obtained for setting a given state feature to zero during a run of Q-GD v1 on the train function. The baseline (none of the features are set to zero) is a final objective of 1.91 and 51/946 half/accept ratio. As a result of zeroing out a state feature, DQN v1 chooses to half the learning rate more frequently and ends up with a worse solution. This experiment shows that DQN v1 depends on each feature to determine the appropriate action. 

 \begin{table}[h] 
 \caption{Effect of setting a state feature to zero. Baseline (none of the feature are set to zero) is a final objective value of 1.91 and a 51/946  half/accept ratio.} \centering 
\begin{tabular}{c c c c} 
\hline\hline
Feature & Objective & Half/Accept \\
\hline
objective value & 1.96 & 320/677 \\
gradient norm & 1.98 & 273/724 \\
alignment measure & 2.04 & 364/633
\end{tabular}
\label{table:features}
\end{table}

\section{Related Work}\label{sec:related}
Neural network models yield state of the art performance in speech recognition, natural language processing, and computer vision \cite{hinton2012deep, krizhevsky2012imagenet,collobert2008unified}. Tesauro popularized neural networks as an approximation to the value function \cite{tesauro1995temporal}, which Riedmiller later extended to the action-value function with the advent of the Neural Fitted Q Iteration \cite{riedmiller2005neural}. Applications of using neural networks in RL appear in settings ranging from playing games to robotics \cite{mnih2015human,lin1993reinforcement}.


Using reinforcement learning to replace an optimization heuristic or be embedded within the optimization algorithm has been explored in a variety of domains \cite{boyan1998learning, dorigo2014ant, miagkikh1999global, moll2000machine,ruvolo2009optimization}. However, none of the previous approaches use deep Q-learning or our proposed RL formulation. Our work is most similar to \cite{ruvolo2009optimization}; the authors use RL to replace a Levenberg-Marquardt heuristic for controlling a damping parameter used in a Gauss-Newton update routine. Unlike our work, they approximate the action-value function by a linear combination of basis functions, which they train using Least Square Policy Iteration. To our knowledge, our work is the first to successfully apply deep Q-learning to controlling an optimization hyperparameter.

\section{Conclusions} \label{sec:conclusions}
This paper lays the foundation for using deep Q-learning to control an optimization hyperparameter. We defined the state, reward function, and actions such that a DQN could learn how to control the learning rate used in a gradient-based optimization routine, resulting in two Q-gradient descent algorithms. Given that there are no theoretical guarantees that the DQN would find the optimal policy or that its q-values would converge, we presented numerical evidence that the Q-GD algorithms performed better than either gradient descent with an Armijo or nonmonotone line search and that the DQNs' q-values for the optimal action converged to the discounted return of rewards at each time step. Additionally, we demonstrated that the Q-GD algorithms were able to generalize when the train function was replaced with a larger test function.

 A main advantage of the Q-gradient descent method is that it can easily incorporate any objective statistic by adding it to the state feature vector. Future areas of work involve using this framework to explore additional state features that can facilitate optimization decisions. We trained the DQNs in a simple environment in order to demonstrate feasibility. To make this method practical for large scale optimization it is necessary to extend Q-GD to the stochastic regime, that is create Q-stochastic gradient descent. A final area of work involves expanding the actions to include controlling additional hyperparameters, such as a momentum term. Overall, the presented framework allows us to develop new optimization algorithms and gain intuition to the type of strategies that are successful for minimizing neural networks.

\newpage
\small
\bibliographystyle{plain}
\bibliography{referencesIBM}

\end{document}